\documentclass{amsart}  
\usepackage{varioref}   
\usepackage[all]{xy}    
\usepackage{graphics}
\usepackage{graphicx}
\newtheorem{thm}{Theorem} 
\newcounter{formelnummer} 
\newcounter{diagramm}     
\newtheorem{cor}[thm]{Corollary} 
\newtheorem{lem}[thm]{Lemma} 
\newtheorem{prop}[thm]{Proposition} 
\theoremstyle{definition} 
\newtheorem{defn}[thm]{Definition} 
\newtheorem{rem}[thm]{Remark} 
\newtheorem{ex}[thm]{Example} 
\newtheorem{notation}[thm]{Notation} 
\numberwithin{equation}{section} 

\newcommand{\captiondiagramm}{\begin{center}\refstepcounter{diagramm}{\it Diagram }\arabic{diagramm}\end{center}}

\newcommand{\s}{\scriptstyle}

\newcommand{\fs}{\footnotesize}

\newcommand{\CC}{\mathbb{C}}
\newcommand{\QQ}{\mathbb{Q}}
\newcommand{\RR}{\mathbb{R}}
\newcommand{\ZZ}{\mathbb{Z}}
\newcommand{\NN}{\mathbb{N}}
\newcommand{\HH}{\mathbb{H}} 
 
\newcommand{\AAA}{\mathbb{A}}
\newcommand{\PP}{\mathbb{P}}

\newcommand{\ww}{\omega} 
\newcommand{\wwdach}{\hat{\omega}}

\newcommand{\Qquer}{\overline{\QQ}}

\newcommand{\absgal}{\mbox{Gal}(\Qquer/\QQ)}

\newcommand{\GT}{\mbox{GT}}

\newcommand{\GTdach}{\widehat{\GT}}

\newcommand{\gal}{\mbox{Gal}}
\newcommand{\emgal}{\mbox{\em Gal}}

\newcommand{\aut}{\mbox{Aut}}
\newcommand{\emaut}{\mbox{\em Aut}}  

\newcommand{\eminaut}{\mbox{\em Inn}}

\newcommand{\id}{\mbox{id}}

\newcommand{\torus}{E^{\star}} 
\newcommand{\toruso}{E} 

\newcommand{\torusb}{E_B^{\star}}
\newcommand{\torusob}{E_B}
\newcommand{\aff}{\mbox{Aff}}
\newcommand{\affplus}{\mbox{Aff}^+}
\newcommand{\emaffplus}{\mbox{\em Aff}^+}
\newcommand{\emaff}{\mbox{\em Aff}}
\newcommand{\veech}{\Gamma}
\newcommand{\hol}{\mbox{\bf hol}}  

\newcommand{\hhol}{\mbox{\bf aff}}  
\newcommand{\fdach}{\hat{f}}

\newcommand{\slzwei}{\mbox{SL}_2}
\newcommand{\emslzwei}{\mbox{\em SL}_2}
\newcommand{\pslzwei}{\mbox{PSL}_2}
\newcommand{\dev}{\mbox{\bf dev}}

\newcommand{\der}{\mbox{\bf der}}

\newcommand{\out}{\mbox{Out}}
\newcommand{\emout}{\mbox{\em Out}}

\newcommand{\stern}{\star}

\newcommand{\fdachstern}{\hat{f}_{\star}}

\newcommand{\proj}{\mbox{proj}}
\newcommand{\emproj}{\mbox{\em  proj}}
\newcommand{\gen}{\mbox{\bf Gen}}
\newcommand{\emgen}{\mbox{\em \bf Gen}}
\newcommand{\rep}{\mbox{\bf Rep}}
\newcommand{\emrep}{\mbox{\em \bf Rep}}
\newcommand{\veechquer}{\bar{\Gamma}}

\begin{document}
\title{An algorithm for finding the Veech group of an origami}
\author{G. Schmith\"usen}
\address{Mathematisches Institut II, Englerstr. 2, 76128 Karlsruhe}
\email{schmithuesen@math.uni-karlsruhe.de}

\subjclass[2000]{14H10, 14H30, 53C10}
\keywords{Teichm\"uller curves, Veech groups, Origami}

\date{\today}

\begin{abstract}
We study the Veech group of an origami, i.e.\ of a translation
surface, tessellated by parallelograms. We show that it is
isomorphic to the image of a certain subgroup of $\aut^+(F_2)$ in
$\slzwei(\ZZ) \cong \out^+(F_2)$. Based on this we present an
algorithm that determines the Veech group.
\end{abstract}
\maketitle

\section{Origamis as Teichm\"uller curves}\label{intro}

(Oriented) origamis (as defined in section \ref{settings}) can
be described as follows: Take finitely many copies of the unit
square in $\CC$ and glue them together such that left edges are
glued with right edges and upper edges with lower ones (compare \cite{l}, \cite{mm}).
This defines a compact surface $S$. We restrict ourselves to the cases where
$S$ is connected. \\
Lifting the structure of $\CC$
via the squares defines a translation structure on $S^{\star} :=
S - \{P_1,\ldots,P_n\}$,
where $P_1,\ldots,P_n$ are  finitely many points on $S$.
One can vary the structure on $S^{\star}$ as follows:
For each $\tau\in \HH$ identify the squares on $S$ with the
parallelogram $P(\tau)$ in $\CC$ defined by the vertices $0,1,\tau,1+\tau$.
This defines an isometric
embedding of the upper half plane $\HH$ into the Teichm\"uller
space $T_{g,n}$, where $g$ is the genus of $S$. This embedding is
described in detail in \cite{l} and \cite{m} in the more general context of Teichm\"uller curves.
The image of $\HH$
in $T_{g,n}$ under this embedding is a complex geodesic $\Delta \subset T_{g,n}$. The image $C$
of $\Delta$ in the moduli space $M_{g,n}$ under the natural projection
$T_{g,n} \rightarrow M_{g,n}$ is birational to the mirror image of
$\HH/{\veech}(O)$ (\cite{l}, \cite{m}), where $\veech(O)$ is the
Veech group of an origami $O$, defined as in section
\ref{settings}.\\
$\HH/\Gamma(O)$ is an algebraic curve defined over
$\bar{\QQ}$ (see \ref{geometry}).
One has even more: The embedded curve $C$ in $M_{g,n}$ is an irreducible component of
a Hurwitz space and thus also defined over $\Qquer$ (\cite{mm}).
In \cite{l}, where origamis were originally
introduced, Pierre Lochak  suggests to study them in the
context of the action of $\absgal$ on combinatorial objects, in some
sense as generalization of the study of dessins d'enfants. The group
$\absgal$  acts on the set of {\it origami curves} in $M_{g,n}$, and
this action is faithful as shown
in \cite{mm}.\\
Origami curves represent a special kind of (imprimitive)
Teichm\"uller curves, described in
\cite{m}, namely those that arise via a torus.\\

In this article we study the Veech group $\veech(O)$ of origamis $O$. We describe an algorithm
that finds generators and coset representatives of $\veech(O)$ in $\slzwei(\ZZ)$ and
calculates the genus and the number of points at infinity of
$\HH/\veech(O)$.

{\bf Acknowledgements:}\\
I would like to thank Pierre Lochak for having me introduced to
this subject and Frank Herrlich, my supervisor, for many helpful
discussions. I am grateful to them and to Markus Even and Martin
M\"oller for valuable remarks and suggestions. Furthermore I thank
Stefan K\"uhnlein for the idea of the proof for
Proposition \ref{cong}.

\section{Veech groups of origamis}\label{basicidea}

The algorithm we want to present is based on the following
Proposition \ref{mainprop}. We denote by $F_2$ the free group in
two generators and by $\aut^+(F_2)$ the group of orientation preserving automorphisms of $F_2$. Furthermore,
we use the fact that $\slzwei(\ZZ)$ is isomorphic to
$\out^+(F_2)$, the group of outer orientation preserving automorphisms 
of $F_2$, and denote by $\hat{\beta}: \aut^+(F_2) \rightarrow \out^+(F_2) \cong \slzwei(\ZZ)$
the canonical projection (see Lemma \ref{lemma}).\\
To an origami $O:=(p:X \rightarrow \torus)$ we will associate a
subgroup $H\cong \gal(\HH/X)$ of $F_2$. $\Gamma(O)$ is the
Veech group of $O$.

\begin{prop}\label{mainprop}
Let $O$ be an origami. Let $\emaff^+(H) := \{\gamma\in \emaut^{+}(F_2)| \gamma(H) = H\}$. Then we have:
\[\veech(O) = \hat{\beta}(\emaff^+(H)) \subseteq \emslzwei(\ZZ).\]
\end{prop}

The aim of Section \ref{basicidea} is to explain the notations and prove the statement 
of Proposition \ref{mainprop}.

\subsection{Origamis, translation surfaces and the Veech Group}\label{settings}\hspace*{\fill}\\

In the following let $\toruso$ be a fixed torus and $\torus :=
\toruso - \{\bar{P}\} \quad (\mbox{for some } \bar{P} \in
\toruso)$ be a once punctured torus.

\begin{defn}
An (oriented) {\em origami} $O$ (of genus $g \geq 1$) is a (topological) unramified covering
$p: X \rightarrow \torus$, where $X$ is obtained by erasing finitely many
points of a compact surface $\bar{X}$ of genus $g$.
\end{defn}

Fix a (topological) unramified universal covering
$u:\tilde{X} \rightarrow X$ of $X$. Then $v := p\circ u$ is a universal covering of $\torus$.\\
Let $\gal(\tilde{X}/\torus)$ be the group of its deck
transformations. It is naturally isomorphic to the fundamental
group $\pi_1(\torus,\bar{Q})$ of $\torus$ with an arbitrary base
point $\bar{Q} \in \torus$. Furthermore, $\pi_1(\torus,\bar{Q})$
is isomorphic to $F_2 := F_2(x,y)$, the free group in the two
generators $x$ and $y$. Fix this isomorphism $\alpha: F_2
\rightarrow \pi_1(\torus,\bar{Q}) \stackrel{\mbox{\scriptsize
can}}{\cong} \gal(\tilde{X}/\torus)$ such that $\alpha(x)$ and
$\alpha(y)$
define  a canonical marking on $\torus$.\\
Then, $H :=
\gal(\tilde{X}/X) \subseteq \gal(\tilde{X}/\torus)$  is considered (via
(can$\, \circ\, \alpha)^{-1}$) as subgroup of $F_2$.
\begin{notation}\label{H}
\[H := \gal(\tilde{X}/X) \subseteq \gal(\tilde{X}/\torus) = F_2(x,y) =: F_2.\]
\end{notation}

We will consider translation structures on $X$ induced by translation structures on $\torus$.
Therefore we first want to recall some definitions and notations (see \cite{th}, \cite{gj}).\\

An atlas on a surface $X$ such that all transition maps are translations defines a {\em translation structure}
$\mu$ on $X$. $X_{\mu} := (X,\mu)$ is called {\em translation surface}. We call
\[\affplus(X_{\mu}) := \{f:X_{\mu} \rightarrow X_{\mu}|
f \mbox{ is an orientation preserving affine
diffeomorphism\footnotemark} \}\]
the affine group of $X_{\mu}$.
\footnotetext{In the following all diffeomorphisms are orientation preserving}\\
Let $u: \tilde{X} \rightarrow X$ be a (topological)
universal covering of $X$. Then $\tilde{X}$ becomes a translation surface $\tilde{X}_{\eta}$
by lifting the structure $\mu$ on $X$ via $u$ to $\eta$ on
$\tilde{X}$. A fixed chart $(U,\eta_U)$ of
$\tilde{X}_{\eta}$ defines a holomorphic map $\dev: \tilde{X}_{\eta} \rightarrow
\CC$ ({\em developing map}) such that
\[ \eta_U = dev|_{U} \quad \mbox{ and } \quad
 \eta_{U'} = t \circ dev|_{U'}   \mbox{ for a translation }  t :=  t(U',\eta_{U'})\]
for any other chart $ (U',\eta_{U'})$ of $\tilde{X}_{\eta}$.\\
For any affine diffeomorphism $\hat{f}$ of $\tilde{X}_{\eta}$ there is a unique affine diffeomorphism
$\hhol(\fdach)$ of $\CC$ such that $\dev \circ \fdach = \hhol(\fdach) \circ \dev$. We call $\hhol$ the group homomorphism
\[\hhol : \affplus(\tilde{X}_{\eta}) \rightarrow  \affplus(\CC), \fdach \mapsto \hhol(\fdach).\]
The {\em holonomy mapping} $\hol$ is the restriction of $\hhol$ to the subgroup $H = {\gal(\tilde{X}/X)}$ of
$\aff^+(\tilde{X})$. If $\proj$ is the natural projection $\proj:\affplus(\CC) \rightarrow \mbox{GL}_2(\RR)$, then
the group homomorphism
\begin{eqnarray*}
\der : \affplus(X_{\mu}) \rightarrow \mbox{GL}_2(\RR), f \mapsto \proj(\hhol(\fdach))\\
\mbox{where $\fdach$ is some lift of $f$ to $\tilde{X}$}
\end{eqnarray*}
is well defined and called {\em derived map} .\\

$\veech(X_{\mu}) := \der(\affplus(X_{\mu})) \subseteq \mbox{GL}_2(\RR)$ is called  the {\em Veech group} of $X_{\mu}$.
It is independent of the choice of the chart $(U,\mu_U)$ which we used to define $\dev$. If $X$ is precompact, i.e.
$X$ is obtained by erasing finitely many points from a compact Riemann surface $\bar{X}$,
then every
$f \in \affplus(X_{\mu})$ preserves the volume. Thus, $\Gamma(X_{\mu})$ is in $\slzwei(\RR)$. \\

Now, given an origami $O = (p:X\rightarrow \torus)$ as above, any
matrix
\[B = \begin{pmatrix}a&b\\c&d \end{pmatrix} \in \slzwei(\RR)\]
defines a translation structure on $X$ as follows:\\
Take the lattice
\[\Lambda_B :=
    <\vec{v}_1 := \begin{pmatrix}a\\c\end{pmatrix},
     \vec{v}_2 := \begin{pmatrix}b\\d\end{pmatrix}>
   \mbox{ in } \CC.\]
Let $\torusob := \CC/{\Lambda_B}$ be the elliptic curve defined by
$\Lambda_B$  and let $\torusb$ be the once punctured elliptic
curve (obtained by erasing the image of $0$ from $\torusob$) with
the induced translation structure. Fix some point $Q_B$ in $\CC -
\Lambda_B$. Let $\bar{Q}_B$ be its image on $\torusb$. Fix
furthermore as canonical marking the images of the segments from
$Q_B$ to $Q_B + \vec{v}_1$ and from $Q_B$ to $Q_B + \vec{v}_2$ on
$\torusb$. Identify $\torusb$ with $\torus$ via a diffeomorphism
respecting the canonical markings. This way $p$ defines an
unramified covering of $\torusb$. Let $\mu_B$ be the translation
structure on $X$ defined by lifting the translation structure on
$\torusb$ to $X$ via $p$ ($\mu_B$ depends also on $p$ !).
Similarly let $\eta_B$ be the translation structure on the fixed
universal covering $\tilde{X}$ defined via $u$.

\begin{notation}
Denote by $X_B := X_B(O) :=  (X,\mu_B)$ the surface $X$ with translation structure $\mu_B$.
Furthermore, denote by $\tilde{X}_B$ the translation surface $(\tilde{X},\eta_B)$.
\end{notation}

Then the maps $p_B: X_B \rightarrow \torusb$, $u_B: \tilde{X}_B \rightarrow X_B$ and
$v_B: \tilde{X}_B \rightarrow \torusb$ induced by $p$, $u$ and $v$ are translation maps.\\

Let $\dev_B:\tilde{X}_B \rightarrow \CC$ be a developing map of $\tilde{X}_B$ (and thus also for $X_B$ and $\torusb$)
 and $\der_B: \affplus(\tilde{X}_B) \rightarrow \mbox{GL}_2(\RR)$ the corresponding derived map.\\
The proof  of the following Remark \ref{affindependent} shows that the affine group of an origami surface $X_B$ does not depend (up to conjugacy) on the choice of the matrix $B$.

\begin{rem}\label{affindependent}
Let $B, B'$ be in $\slzwei(\RR)$. Then
\[\affplus(X_{B}(O)) \cong \affplus(X_{B'}(O)) \mbox{ and } \Gamma(X_{B'}(O)) = B'B^{-1}\Gamma(X_B(O))BB'^{-1}.\]
\end{rem}

\begin{proof}
The map $\varphi: X_B(O) \rightarrow X_{B'}(O)$ that is topologically the identity on $X$ is an affine diffeomorphism and induces the group isomorphism: \[\affplus(X_B(O)) \rightarrow \affplus(X_{B'}(O)), f \mapsto \varphi\circ f\circ\varphi^{-1}.\]
Since $\der(\varphi) = B'B^{-1}$, we have $\der(\varphi f\varphi^{-1}) = B'B^{-1}\der(f)BB'^{-1}$
\end{proof}

Since the Veech group depends only up to conjugacy on the choice
of $B$, we will restrict to the case of $B = I$, the identity
matrix. If not stated otherwise, we will denote $\tilde{X} := \tilde{X}_I$,
$\der := \der_I$, $\dev := \dev_I$, $X := X_I$, $\toruso :=
\toruso_I$, $\Lambda := \Lambda_I$$, \torus :=
\torus_I$, $\mu := \mu_I$ and $\veech(O) := \veech(X_I(O))$.\\
By the uniformization theorem there exists a biholomorphic map
$\delta: \HH \rightarrow \tilde{X} = \tilde{X}_I$, where $\HH$ is the complex upper
half plane. $\HH$ becomes via $\delta$ a
translation surface. We will identify in the following $\HH$ with
$\tilde{X} = \tilde{X}_I$.\\

\begin{prop}\label{affinegroupoftorus}
Let $O = (p:X\rightarrow \torus)$ be an origami and $\HH$ be the upper half plane,
endowed with the translation structure induced by $O$ as above. Then we have:
\begin{enumerate}
\item \label{subgroup} $\veech(O)$ is a subgroup of $\veech(\HH)$. 
\item
\label{veechfull} $\Gamma(\torus) = \Gamma(\HH) = \emslzwei(\ZZ).$
\item \label{partD} Let $f$ be in $\emaffplus(X)$. Then $f$
descends via $p$ to some $\bar{f} \in \emaffplus(\torus)$ and
Diagram \ref{diagramm1} is commutative with $A := \der(f)$, with
$\fdach$ some lift of $f$ to $\HH$ and with some $b \in \ZZ^2$.
\end{enumerate}

\[ \xymatrix{
   &\HH\ar[rrr]^{\fdach}\ar[dl]_{h} \ar[ddr]^<<<<<<<{u}\ar[d]_{\dev} \ar @{} [drrr] |{(A)} &&
   &\HH\ar[ddl]_<<<<<<<{u}\ar[d]^{\dev}\ar[dr]^{h}&\\
   \CC-\Lambda  \ar@{}[r]|{\subset}  \ar[ddrr]^{w}  &
   \CC\ar[rrr]^{z \mapsto Az + b} \ar @{} [drrr] |{}&&&
   \CC \ar @{} [r] |{\supset}  & \CC-\Lambda  \ar[ddll]_{w} \\
     && X \ar[r]^{f} \ar[d]_{p} \ar@{} [dr] |{}& X \ar[d]^{p}&&\\
     && \torus \ar[r]^{\bar{f}} & \torus &&}
\]
\nopagebreak
\begin{center}\refstepcounter{diagramm}{\it Diagram }\arabic{diagramm}\label{diagramm1}
\end{center}

\end{prop}

\begin{proof}\hspace*{\fill}\\
\underline{\bf 1.:} Let $f$ be in $\affplus(X)$ and $\fdach$ be
some lift of $f$ via $u$. Since the translation structure on $\HH$
is lifted via $u$, $\fdach$ is also affine and $\der(\fdach) =
\der(f)$. Hence, $\Gamma(O) \subseteq \Gamma(\HH)$.\\[.3cm]
\underline{\bf 2.:} Let $\CC \rightarrow \toruso$ be the universal
covering and $w: \CC - \Lambda \rightarrow \torus$ its restriction
to $\CC - \Lambda$. Since $v = p\circ u$ is the universal covering of
$\torus$, there is an unramified covering $h:\HH \rightarrow \CC -
\Lambda$, such that $w\circ h = v = p\circ u$. But since the
structure on $\HH$ was obtained by lifting the translation
structure on $\torus$ via $v$, this map $h$ is locally a chart of
$\HH = \tilde{X}_I$. Thus, the map $h$ is a developing map and the image of
this developing map $\dev$ is $\CC - \Lambda$.\\
Now, let $A$ be in $\Gamma(\HH)$, hence $A = \der(\fdach)$ for
some $\fdach \in \affplus(\HH)$. By the definition of $\der$ and
$\dev$ Part (A) of Diagram \ref{diagramm1} is commutative for some
$b \in \ZZ^2$, i.\ e.
\[ (z \mapsto Az + b)\circ \dev = \dev \circ \fdach.\]
Since the image of $\dev$ is in $\CC - \Lambda$, the map $z
\mapsto Az + b$ respects $\Lambda = \ZZ^2$. Thus, $A$ is in
$\slzwei(\ZZ)$.
Hence, we have: $\Gamma(\HH) \subset \slzwei(\ZZ)$.\\
Conversely, taking a matrix $A \in \slzwei(\ZZ)$ the map $z
\mapsto Az$  descends to an affine diffeomorphism $\bar{f} \in
\aff^+(\torus)$. This can be lifted to some $\fdach \in
\aff^+(\HH)$ with $\der(\fdach) = A$. Thus, we have: $\slzwei(\ZZ)
\subset \Gamma(\HH)$.\\
Using the same arguments it follows that also $\Gamma(\torus) =
\slzwei(\ZZ)$.\\[.3cm]
\underline{\bf 3.:} Let $\fdach \in \aff^+(\HH)$ be some lift of
$f$ to $\HH$. By the proof of (2) it follows that $\fdach$
descends via $w\circ h = v$ to some $\bar{f} \in \aff^+(\torus)$
and that Diagram \ref{diagramm1} is commutative.
\end{proof}

From (\ref{subgroup}) and (\ref{veechfull}) of Proposition \ref{affinegroupoftorus} we see in particular that 
the Veech group 
$\veech(O)$ of an origami $O$ is always a subgroup of $\slzwei(\ZZ)$. It follows
from \cite[Thm.\ 5.5]{gj}, that it has finite index in $\slzwei(\ZZ)$. This result will play 
a crucial role in section \ref{final}.\\

An immediate consequence of Proposition \ref{affinegroupoftorus} is

\begin{cor}\label{zwischenergebnis} 
\[\Gamma(O) = \{A \in \emslzwei(\ZZ)| A = \der(\fdach) \mbox{ \em for some } \fdach \in \emaffplus(\HH) \mbox{\em \ that
 descends to $X$ via } u  \}\\. \label{gammao}
\]
\end{cor}

To prove Proposition
\ref{mainprop} from Corollary \ref{gammao} we have to state a condition for $\fdach$  in
$\affplus(\HH)$ to descend via $u$ to some $f \in \affplus(X)$.

\subsection{When does an element in $\boldmath{\affplus(\HH)}$ descend to $X$?}\label{descend}\hspace{\fill}\\

Recall that $H = \gal(\HH/X)\subset F_2 = \gal(\HH/\torus)
\subseteq \pslzwei(\RR)$ (Notation \ref{H}). We define the group
homomorphism
\begin{eqnarray*}
\stern \, : \quad \aff^+(\HH) &\rightarrow& \aut^+(F_2)\\
    \fdach &\mapsto& (\fdachstern: \sigma \mapsto \fdach\circ \sigma\circ
    \fdach^{-1})
 \end{eqnarray*}
Remark that
\[ F_2 = \gal(\HH/\torus) = \{\fdach\in \aff^+(\HH)| \der(\fdach) = I\}.
\hspace*{2cm }\refstepcounter{formelnummer}[\arabic{formelnummer}]
\label{kern}
\]
 The map $\star$ \; is well defined, since
$\fdach\circ\sigma\circ\fdach^{-1}$ is again affine with the
derivative $\der(\fdach)\cdot I\cdot\der(\fdach)^{-1} = I$ and thus in $F_2$.
\begin{lem} \label{lemma}
We have the following properties of $\stern$\;:\nopagebreak
\begin{enumerate}\nopagebreak
\item  \label{eins} The following two sequences are exact and the
diagram is commutative:
\[ \xymatrix{
   1 \ar[r]& F_2 \ar[r]\ar[d]^{\alpha}_{\cong} \ar @{} [dr] |{(A)}&
   \emaff^+(\HH) \ar[r]^{\der} \ar[d]^{\stern}_{\cong} \ar @{} [dr] |{(B)}&
   \emslzwei(\ZZ)  \ar[r] & 1\\
   1 \ar[r]& {\mbox{\eminaut}(F_2)} \ar[r]&
   \emaut^+(F_2) \ar[r] &
   \emout^+(F_2) \ar[u]^{\cong}_{\beta} \ar[r] & 1}
\]
\captiondiagramm Here, $\mbox{\eminaut}(F_2)$ is the group of inner
automorphisms
of $F_2$, $\alpha$ is the natural isomorphism $F_2
\rightarrow \mbox{\eminaut}(F_2), x \mapsto (y \mapsto xyx^{-1})$ ,
$\beta: \emout^{+}(F_2) \rightarrow \emslzwei(\ZZ)$ is the group
isomorphism induced by the natural homomorphism:
\[\hat{\beta}:\emaut^+(F_2) \rightarrow \mbox{\em SL}_2(\ZZ), \,
\varphi \mapsto A := \begin{pmatrix} a&b\\c&d\end{pmatrix},\]
where $a$ is the number of $x$ appearing in $\varphi(x)$, $b$ the
number of $x$ appearing in $\varphi(y)$, $c$ the number of $y$ in
$\varphi(x)$ and $d$ the number of $y$ in $\varphi(y)$ (see
\cite[I 4.5, p.25]{ls}). Recall that for the canonical
projection $\emproj: F_2 \rightarrow \ZZ^2$ sending $x$ to
$(1,0)^t$ and $y$ to $(0,1)^t$ one has:
\[ \hspace*{2cm}
\forall \varphi \in \emaut^+(F_2), A := \hat{\beta}(\varphi)
\quad\quad \emproj\,\circ\, \varphi = (z \mapsto A\cdot z)\, \circ \, \emproj.
\hspace*{2cm}\refstepcounter{formelnummer}[\arabic{formelnummer}]
\label{gleichung1}
\]
\item
An element $\fdach \in \emaff^+(\HH)$ descends to $X$ via $p$ iff $\fdach_{\stern}(H) = H$.
\end{enumerate}
\end{lem}

\begin{proof}\hspace*{\fill}\\
{\bf \underline{1.}}:\\[.2cm]
The exactness of the first sequence follows by Equation \ref{kern}
and by Proposition \ref{affinegroupoftorus}. The exactness of the
second sequence is true by the definition of $\out^{+}(F_2)$.\\
 The commutativity of Part $(A)$ of the Diagram is true
by definition of $\star$. We prove now the commutativity of Part
$(B)$:\\
We have chosen the isomorphism $F_2 = F_2(x,y) \cong
\gal(\HH/\torus)$ and the translation structure on $\torus =
\torus_I$ in such a way that:
\[ \hhol(x) = (z \mapsto z + \begin{pmatrix}1\\0 \end{pmatrix}) \mbox{ and }
\hhol(y) = (z \mapsto z + \begin{pmatrix}0\\1 \end{pmatrix}).\]
Thus, $\hhol|_{F_2} (= \hol)$ is the natural projection $\proj: F_2 \rightarrow \ZZ^2$. Here we identify
the group of translations of $\CC$ along some vector in $\ZZ^2$ canonically with $\ZZ^2$.\\
Consider the following diagram:
\[ \xymatrix{F_2 \ar[r]^{\fdachstern} \ar[d]_{\mbox{proj}}&
         F_2 \ar[d]^{\mbox{proj}}\\
       \ZZ^2 \ar[r]^{z \mapsto A\cdot
         z} & \ZZ^2  }
\]
\begin{center}
\refstepcounter{diagramm}{\it Diagram }\arabic{diagramm}\label{diagramm2}
\end{center}

Diagram \ref{diagramm2} is commutative with $A := \der(\fdach)$:\\
Let $\sigma$ be in $F_2 = \gal(\HH/\torus)$. We have to show that
$\mbox{proj}(\fdachstern(\sigma)) = A\cdot\mbox{proj}(\sigma)$.\\
We have $\hhol(\sigma) = (z \mapsto z + c)$ and  $\hhol(\fdach) = (z \mapsto Az+b)$ for some
$b, c \in \ZZ^2$.
Thus we get:
\[\proj(\fdachstern(\sigma)) = \hhol(\fdachstern(\sigma))  =
     \hhol(\fdach)\hhol(\sigma)\hhol(\fdach^{-1}) = (z  \mapsto z + Ac).
\]
Hence, Diagram \ref{diagramm2} is commutative with $A = \der(\fdach)$.\\
To conclude we use  that Diagram \ref{diagramm2} is also commutative with $A = \hat{\beta}(\fdachstern)$
(see equation [\ref{gleichung1}]). Thus, $\der(\fdach) = \hat{\beta}(\fdachstern)$ and ({\em B}\,) is
commutative.\\
Finally, $\alpha$ and $\beta$ are both isomorphisms, thus $\star$
is also an isomorphism.\\ [.3cm]

{\bf \underline{2.}}:
$\fdach \in \affplus(\HH)$ descends to $X$ via $p$ $\Leftrightarrow$ for all $z \in \HH, \sigma \in H = \gal(\HH/X)$
there is some $\tilde{\sigma}_{z,\sigma} \in H$ such that $\tilde{\sigma}_{z,\sigma}(\fdach(z)) = \fdach(\sigma(z))$.\\
For $\tilde{\sigma} := \fdachstern(\sigma)$ we have by definition of $\fdachstern$:
 $\tilde{\sigma}(\fdach(z)) = \fdach(\sigma(z))$ for all $z \in \HH$.
Since $F_2$ operates fixpointfree on $\HH$ it follows from the
last equation that $\tilde{\sigma}_{z,\sigma}$ has to be equal to
$\tilde{\sigma} = \fdachstern(\sigma)$ . On the other hand,
$\tilde{\sigma}_{z,\sigma}$ has to be in  $H$. This proves (2).
\end{proof}

Now Proposition \ref{mainprop} follows from  Corollary \ref{gammao} and Lemma \ref{lemma}. 
\hspace*{\fill}$ _{\Box}$\\

As result of Proposition \ref{mainprop} we get: In order to check whether $A \in \slzwei(\ZZ)$
is in $\Gamma(O)$, we
have to check if there exists a lift $\gamma_A \in \aut^+(F_2)$ of $A$ (i.e.\ 
a preimage of $A$ under
$\hat{\beta}$) that fixes $H$. The following Corollary translates
this into a finite problem that can be left to a computer.

\begin{cor} (to Proposition \ref{mainprop})\label{cor}\\
Let $O = (p: X \rightarrow \torus)$ be an origami of degree $d$, $F_2 = \emgal(\HH/\torus)$,
$H = \emgal(\HH/X)$ as above. Let $h_1, \ldots, h_k$ be generators of $H$ and
${\sigma}_1, \ldots, {\sigma}_d $ a system of  right coset representatives of
$H\backslash F_2$ (denote the right coset $H\cdot{\sigma}_i$ by $\bar{\sigma}_i$).\\
Further let $\gamma_A^0 \in \emaut^+(F_2)$ be some fixed lift of $A \in \emslzwei(\ZZ)$. Then
\[A \in \Gamma(O) \Leftrightarrow \exists \, i \in \{1,\ldots, d\} \mbox{ such that }
\bar{\sigma}_i\cdot\gamma^0_A(h_j) = \bar{\sigma}_i \mbox{ for all } j \in \{1,\ldots , k\}.\]
\end{cor}

\begin{proof}
Let $\gamma_A$ be another lift of $A$. Thus $\gamma^0_A = 
\sigma^{-1}\cdot\gamma_A\cdot\sigma$ for some $\sigma \in F_2$ and we have for all $h$ in $H$:
\[\gamma_A(h) \in H \Leftrightarrow \sigma\cdot\gamma^0_A(h)\cdot\sigma^{-1}\in H \Leftrightarrow H\cdot\sigma\cdot\gamma^0_A(h) = 
H\cdot\sigma \Leftrightarrow \bar{\sigma}\cdot\gamma^0_A(h)=\bar{\sigma}\]
Hence, the claim follows from Proposition \ref{mainprop}.
\end{proof}

\section{The algorithm}

Let $O = (p:X\rightarrow \torus)$ be a given origami of degree $d$.
In this section we present our algorithm that
determines the Veech group $\Gamma(O)$.
We have subdivided this description into four parts: In \ref{lift} we describe how to find some
lift $\gamma_A \in \aut^+(F_2)$ for any matrix
$A$ in $\slzwei(\ZZ) \cong \out^+(F_2)$, in \ref{descide} we show how
to decide whether a given matrix
$A \in \slzwei(\ZZ)$ is in $\Gamma(O)$, in \ref{final} we give an algorithm
that determines generators and a system
of coset representatives of $\Gamma(O)$ in $\slzwei(\ZZ)$, and finally in \ref{geometry} we state how to calculate
the genus and the points at infinity of the corresponding {\it Veech curve} $\HH/\Gamma(O)$.\\
In order to illustrate the algorithm we will use the example $O = L(2,3)$.

\begin{minipage}{\linewidth}
\begin{ex}(The Origami $O = L(2,3)$)\label{l23}\\[.3cm] 
\hspace*{2cm}
\begin{minipage}[c]{4cm}

\setlength{\unitlength}{1cm}
\begin{picture}(4,2.3)

\put(-1.5,1.2){$X =$}
\put(0,0){\framebox(1,1){$\star$}}
\put(1,0){\framebox(1,1){}}
\put(2,0){\framebox(1,1){}}
\put(0,1){\framebox(1,1){}}

\thicklines 
\multiput(0,0.5)(0.3,0){10}{\line(1,0){0.2}}
\multiput(0,1.5)(0.3,0){3}{\line(1,0){0.2}} \put(0.9,1.5){\line(1,0){0.1}}

\multiput(0.5,0)(0,0.3){7}{\line(0,1){0.2}} 
\multiput(1.5,0)(0,0.3){3}{\line(0,1){0.2}} \put(1.5,0.9){\line(0,1){0.1}}
\multiput(2.5,0)(0,0.3){3}{\line(0,1){0.2}} \put(2.5,0.9){\line(0,1){0.1}}
\put(1.25,0.5){\vector(1,0){0.3}}
\put(0.5,1.25){\vector(0,1){0.3}}

\put(.5,-0.3){$a$}
\put(1.5,-0.3){$b$}
\put(2.5,-0.3){$c$}
\put(3.1,0.5){$d$}
\put(2.5,1.1){$c$}
\put(1.5,1.1){$b$}
\put(1.1,1.5){$e$}
\put(0.5,2.1){$a$}
\put(-0.3,1.5){$e$}
\put(-0.3,0.5){$d$}

\put(-0.07,-0.07){$\bullet$}
\put(0.93,-0.07){$\bullet$}
\put(1.9,-0.07){$\Diamond$}
\put(2.93,-0.07){$\bullet$}
\put(2.93,0.93){$\bullet$}
\put(1.9,0.93){$\Diamond$}
\put(0.93,0.93){$\bullet$}
\put(0.93,1.93){$\bullet$}
\put(-0.07,1.93){$\bullet$}
\put(-0.07,0.93){$\bullet$}

\put(0.25,0.55){$\sf Q$}

\put(3.5,0){\vector(3,-1){1.5}}
\put(4.5,-0.17){$p$}
\put(5.5,-01.2){$\torus = $}
\put(7,-1.5){\framebox(1,1){}}
\put(6.7,-1){$a$}
\put(7.5,-1.8){$b$}
\put(8.1,-1){$a$}
\put(7.5,-0.4){$b$}
\multiput(7,-1)(0.3,0){3}{\line(1,0){0.2}} \put(7.9,-1){\line(1,0){0.1}}
\multiput(7.5,-1.5)(0,0.3){3}{\line(0,1){0.2}} \put(7.5,-0.6){\line(0,1){0.1}}
\put(7.8,-1){\vector(1,0){0.3}}
\put(7.65,-0.99){$\fs x$}
\put(7.5,-0.75){\vector(0,1){0.3}}
\put(7.6,-0.7){$\fs y$}
\put(7.15,-1.33){$\fs \bar{Q}$}
\put(7.42,-1.05){$\star$}
\put(6.7,-1.85){$\bar{P}$}

\end{picture}
\end{minipage}\\
\vspace*{1.5cm}
\begin{center}
\captiondiagramm \label{l23bild}
\end{center}
\vspace*{.3cm}
\end{ex}
\end{minipage}

In Example \ref{l23} the edges labelled with the same letters are
glued together. This way $X$ becomes a surface of genus 2. The
squares describe the covering map to $\torus$. The point $\bar{P}
\in \toruso$ (at infinity) has 2 preimages on the surface $X$ (the
points $\bullet$ and $\Diamond$),
the degree $d$ of $p$ is 4.\\
We identify $F_2 = \gal(\HH/\torus)$ with the fundamental group of $\torus$
(with base point $\bar{Q}$) and $H = \gal(\HH/X)$ with the fundamental group of $X$
(with base point $Q$).
The projection of the closed paths on $X$ to $\torus$ defines the embedding of $H$ into $F_2$,
$x$ and $y$ are the fixed generators of $F_2$ on $\torus$.
Since the $L(2,3)$-shape is simply connected, the generators of $H$ are obtained by the identifications of the edges.
Thus, $H = <x^3, x^2yx^{-2}, xyx^{-1},yxy^{-1},y^2>$. The index $[F_2 : H]$ is equal to $d = 4$.

\subsection{Lifts from $\boldmath{\slzwei(\ZZ)}$ to the automorphism group of $\boldmath{F_2}$}\label{lift}\hspace*{\fill}\\

Let
\[S := \begin{pmatrix} 0&-1\\1&0\end{pmatrix} \mbox{ and }
  T := \begin{pmatrix} 1& 1\\0&1\end{pmatrix}.\]
We will use the fact that $\slzwei(\ZZ)$ is generated by $S$ and $T$ and that $S^{-1} = -S$ and $T^{-1} = -STSTS$. Thus, every $A \in \slzwei(\ZZ)$ can be written as $A = W(S,T)$ or $A = -W(S,T)$, where $W$ is a word in the letters $S$ and $T$.\\
The homomorphisms
\begin{eqnarray*}
\gamma_S&:& F_2 \rightarrow F_2 \mbox{ defined by } \gamma_S(x) = y
   \mbox{ and } \gamma_S(y) = x^{-1} \mbox{, }\\
\gamma_T&:& F_2 \rightarrow F_2 \mbox{ defined by } \gamma_T(x) = x
   \mbox{ and } \gamma_T(y) = xy \quad \mbox{ and }\\
\gamma_{-I}&:& F_2 \rightarrow F_2 \mbox{ defined by } \gamma_{-I}(x) = x^{-1}
   \mbox{ and } \gamma_{-I}(y) = y^{-1}\\
\end{eqnarray*}
are in $\aut^+(F_2)$ with $\hat{\beta}(\gamma_S) = S$, $\hat{\beta}(\gamma_T)=T$
and $\hat{\beta}(\gamma_{-I})=-I$, where the morphism
$\hat{\beta}: \aut^+(F_2) \rightarrow \slzwei(\ZZ)$ is the projection defined in
\ref{descend} (Lemma \ref{lemma}).\\
Hence, for $A = \pm W(S,T)$ the automorphism
$\gamma_A := \pm W(\gamma_S,\gamma_T) \in \aut^+(F_2)$ is a lift of $A$. Hereby we denote
$-W(\gamma_S,\gamma_T) := \gamma_{-I}\circ W(\gamma_S,\gamma_T)$.\\

In order to find a word $W$ such
that $A = W(S,T)$ or $A = -W(S,T)$ we will define a sequence $A_1 := A$, $A_2$, $\ldots$, $A_N$
such that (for $1\leq n < N$) 
\[ A_{n+1} = A_n\cdot T^{-k_n}\cdot S \; (\mbox{with } k_n \in \ZZ)
\mbox{ and } A_N = \pm T^{\pm b_N} \; (\mbox{with } b_N \in \ZZ).\]
From this we get that
$A =\pm T^{\pm b_n}\cdot(-S)\cdot T^{k_{n-1}}\cdot\ldots\cdot(-S)
\cdot T^{k_1}$. We will conclude using that 
$T^{-1} = -STSTS$.\\ 
These considerations give rise to the following 
algorithm, in which we denote 
\[A_n =: \begin{pmatrix}a_n&b_n\\c_n&d_n\end{pmatrix} \mbox{ with } a_n, b_n, c_n, d_n \in \ZZ.\]

\framebox{
\begin{minipage}{\linewidth}
{\bf Algorithm for finding a lift in $\aut^+(F_2)$:} \\

{\it Given:} $A \in \slzwei(\ZZ)$.\\[.2cm]
$n := 1; A_1 := A.$
\begin{enumerate}
\item If $c_n  \neq 0$ find $k_n \in \ZZ$, such that
\[A_{n+1} := A_nT^{-k_n}S \; \mbox{ fulfills } |c_{n+1}| < |c_n| .\]
$k_n := d_n \, \mbox{div} \, c_n$ does this job: 
$d_n = k_nc_n + r_n $ with $r_n \in \{0,1,\ldots,|c_n|-1\}$ 
\[
\Rightarrow A_{n+1} = \begin{pmatrix} -a_nk_n+b_n & -a_n\\r_n & -c_n\end{pmatrix}.\]
Increase $n$ by 1.
\item Iterate Step (1) until $c_n = 0$. Thus \
\begin{eqnarray*}
A_n  &=&
\begin{pmatrix}\pm 1& b_n\\0&\pm 1\end{pmatrix} = \pm T^{\pm b_n} \mbox{ and }\\
A &=& \pm T^{\pm b_n}\cdot(-S)\cdot T^{k_{n-1}}\cdot\ldots\cdot(-S)
\cdot T^{k_1} =:\pm\tilde{W}(S,T,T^{-1}).
\end{eqnarray*}
\item Replace in $\tilde{W}$ each $T^{-1}$ by $-STSTS$\\
  $\Rightarrow$ Word $W$ in $S$ and $T$ with $A = W(S,T)$ or $A = -W(S,T)$.
\item Compute $\gamma_A := W(\gamma_S,\gamma_T)$ or $\gamma_A := -W(\gamma_S,\gamma_T)$.\\
\end{enumerate}
{\em  Result:} $\gamma_A \in \aut(F_2)$ with $\hat{\beta}(\gamma_A) = A$.
\end{minipage}\\
}

\begin{ex} \label{exlift}
\begin{eqnarray*}
\begin{pmatrix}-3&5\\-2&3 \end{pmatrix} =  -T^2ST^3STS
&\Rightarrow&\gamma^0_A = \gamma_{-I}\gamma_{T}^2\gamma_{S}\gamma_{T}^3\gamma_{S}\gamma_{T}\gamma_{S} \\
&\Rightarrow&\gamma^0_A: \; x \mapsto x^{-2}y^{-1}x^{-2}y^{-1}x^{-2}y^{-1}xyx^2 , \quad
              y \mapsto x^{-1}yx^2yx^2yx^2
\end{eqnarray*}
\end{ex}

\subsection{Decide whether $A$ is in the Veech Group $\boldmath{\Gamma(O)}$}\label{descide}\hspace*{\fill}\\

Let $A$ be in $\slzwei(\ZZ)$. We want to decide whether $A$ is in $\Gamma(O)$ or not. As in Corollary \ref{cor} let $h_1,\ldots, h_k$ be generators of $H = \gal(\HH/X) \subseteq F_2 = \gal(\HH/\torus)$, $\sigma_1, \ldots, \sigma_d$ a system of right coset representatives of $H$ in $F_2$ ($\bar{\sigma_i} := H\cdot\sigma_i$) and $\gamma^0_A$ some fixed lift of $A$ in $\aut^+(F_2)$.\\
Corollary \ref{cor} suggests how to build the algorithm:
\[A \in \Gamma(O) \Leftrightarrow \exists i \in \{1,\ldots,d\} \mbox{ such that } \forall j \in \{1,\ldots, k\} \quad \bar{\sigma}_i\cdot \gamma^0_A(h_j) = \bar{\sigma}_i.\]
Hence, the {\it main step} will be to decide for some $\tau \in F_2$ whether
\[\bar{\sigma_i}\cdot\tau = \bar{\sigma_i}.\]

In order to do this we present the origami $O$ as directed graph
$G$ with edges labelled by $x$ and $y$ (see Figure
\ref{l23tree}). The cosets
$\bar{\sigma}_1, \ldots,\bar{\sigma}_d$ are the vertices of $G$.
Each vertex $\bar{\sigma}_i$ is start point of one $x$-edge and
one $y$-edge. The endpoint is $\overline{\sigma_i\cdot x}$ and
$\overline{\sigma_i\cdot y}$, respectively.
\[
\SelectTips{cm}{}
\xymatrix @-1pc {
*++[o][F-]{\bar{y}}\ar@/_/[d]_y \ar@(dr,ur)[]_{x}&&\\
*++[o][F-]{\bar{\mbox{id}}}\ar[u]_y\ar[r]^{x}&
*++[o][F-]{\bar{x}}\ar[r]^{x}\ar@(ur,ul)[]_y&
*++[o][F-]{\bar{x^2}}\ar@(ur,ul)[]_y \ar@/^1pc/[ll]^x
}
\]
\makebox{
\begin{minipage}{\linewidth}
\begin{center}
\refstepcounter{diagramm}
{\it Figure }\arabic{diagramm}: Graph for $O = L(2,3)$.\label{l23tree}
\end{center}
\end{minipage}
}

Writing  $\tau \in F_2$ as word in $x$,$y$,$x^{-1}$ and $y^{-1}$ defines a not necessarily oriented path in $G$ 
starting at the vertex $\bar{\sigma_i}$ with end point $\bar{\sigma}_i\cdot\tau$ . We have:
\[\bar{\sigma_i}\cdot\tau = \bar{\sigma_{i}} \Leftrightarrow \mbox{  this path is closed.}\]

Thus we get the following algorithm.\\

\framebox{
\begin{minipage}{\linewidth}
{\bf Algorithm for deciding whether $\boldmath{A}$ is in $\boldmath{\Gamma(O)}$:}\\

{\em Given:} $A \in \slzwei(\ZZ)$.\\[.2cm]
Calculate some lift $\gamma^0_A \in \aut^+(F_2)$ of $A$ (see \ref{lift}).\\
For $j = 1$ to $k$ do: $\tilde{h}_j := \gamma^0_A(h_j)$.\\
result $:=$ false.\\
for $i=1$ to $d$ do\\[.3cm]
\hspace*{0.5cm}
\begin{minipage}{12cm}
help $:=$ true.\\
for $j = 1$ to $k$ do: if $\bar{\sigma}_i\cdot\tilde{h}_j \neq \bar{\sigma}_i$ (main step, see above) then
help $:=$ false.\\
if help $=$ true then result $:=$ true.\\
\end{minipage}\\
{\em Result:} If the variable 'result' is true, then $A \in \Gamma(O)$, else $A \not\in \Gamma(O)$.
\end{minipage}
}

\begin{ex} (for $O = L(2,3)$)

Let $A := \begin{pmatrix} 1&0\\2&1 \end{pmatrix}.$ Take the lift:
\[ \gamma^0_A: \; x \mapsto xyxyx^{-1} =: u \quad
                  y \mapsto xyxyx^{-1}y^{-1}x^{-1} =: v\]

Generators of $H$ (see Ex.\ \ref{l23}) are:
\[h_1 := x^3, h_2 := xyx^{-1}, h_3 := x^2yx^{-2}, h_4 := yxy^{-1}, h_5 := y^2.\]
For example
$\bar{\id}\cdot \gamma^0_A(h_2) = \bar{\id}\cdot uvu^{-1} = \bar{x}vu^{-1} = \bar{x^2}u^{-1}
= \bar{x^2} \; \Rightarrow \; \gamma_A^0(H) \neq H$.\\
But one has: $\bar{x}\cdot\gamma^0_A(h_i) = \bar{x} \; \forall i\in\{1,\ldots,5\}$.\\
$\Rightarrow \gamma_A(H) = H$ for $\gamma_A = x\cdot\gamma_A^0\cdot x^{-1}$  and $A \in \veech(O)$.
\end{ex}

\subsection{Generators and Coset Representatives of $\boldmath{\veech(O)}$}\label{final}\hspace*{\fill}\\

Let $\veechquer(O)$ be the {\it projective Veech group}, i.e.\ the image of $\Gamma(O)$
under the projection of $\slzwei(\ZZ)$ to $\pslzwei(\ZZ)$. We first give an algorithm that
calculates a list $\gen$ of generators and a list $\rep$ of right coset representatives of
$\veechquer(O)$ in $\pslzwei(\ZZ)$, then we determine $\veech(O)$. The way how we proceed
is based on the Reidemeister-Schreier method (\cite{ls}, II.4).\\
We denote by $\bar{A}$ the image of an element $A\in \slzwei(\ZZ)$ under the projection to
$\pslzwei(\ZZ)$ and, conversely, denote for $\bar{A}$ in $\pslzwei(\ZZ)$ by $A$ some lift of
$\bar{A}$. Moreover, we write $A \sim B$ (respectively $\bar{A} \sim \bar{B}$) if they are in the
same coset, i.e.\ $\veech(O)\cdot A = \veech(O)\cdot B$ (respectively
$\bar{\veech}(O)\cdot \bar{A} = \bar{\veech}(O)\cdot\bar{B}$).\\
Each element of $\pslzwei(\ZZ)$ can be presented as word in
$\bar{S}$ and $\bar{T}$. We use the directed infinite tree shown
in Figure \ref{trees}: The vertices $v_0, v_1, v_2, \ldots$ of the
tree are labelled by elements of $\pslzwei(\ZZ)$. The root $v_0$
is labelled by $\bar{I}$, the image of the identity matrix. Each
vertex is
starting point of two edges, one labelled by $\bar{S}$, one labelled by $\bar{T}$.\\
Each element of $\pslzwei(\ZZ)$ occurs as label of at least one vertex. Starting with $v_0$
we will visit each vertex $v$ (with label $\bar{B}$) and check if it is not yet represented
by the list $\rep$. In this case we will add it to $\rep$. Otherwise for each $\bar{D}$
in $\rep$ that is in the same coset as $\bar{B}$, we add $\bar{B}\cdot\bar{D}^{-1}$ to the
list $\gen$ of generators.

\makebox{
\begin{minipage}{\linewidth}
\[
\SelectTips{cm}{}
\xymatrix @-1pc {
&&&&*++[o][F-]{\scriptstyle\bar{I}}\ar[lld]_{\bar{T}}\ar[rrd]^{\bar{S}}\ar@{}[d]|{v_0} &&&\\
&&*++[o][F-]{\s \bar{T}}\ar[ld]_{\bar{T}}\ar[rd]^{\bar{S}}\ar@{}[d]|{v_1}&&
&&*++[o][F-]{\s \bar{S}}\ar@{.>}[ld]_{\bar{T}}\ar[rd]^{\bar{S}}\ar@{}[d]|{v_2}&\\
&*++[o][F-]{\s \bar{T}^2}\ar@{.>}[ld]_{\ldots}\ar@{}[d]|{v_3}&
&*++[o][F-]{\s \bar{T}\bar{S}}\ar@{}[d]|{v_4}&&
*++[o][F-]{\s \bar{A}_l}\ar[rd]^{\bar{S}}\ar@{}[d]|{v_l}&&
*++[o][F-]{\s \ldots}\\
*++[o][F-]{\s \bar{A}_j}\ar[rd]^{\bar{S}}\ar@{}[d]|{v_j}&\ldots &\ldots &\ldots &\ldots &\ldots
&*++[o][F-]{\s \ldots}\ar[ld]_{\bar{T}}&\ldots &\\
&*++[o][F-]{\s \ldots}\ar[ld]_{\bar{T}}&\ldots &\ldots
&\ldots &*++[o][F-]{\s \bar{A}_m}\ar@{}[d]|{v_m}&\ldots &\ldots &\\
*++[o][F-]{\s \bar{A}_{n+1}}&\ldots &\ldots &\ldots &\ldots &\ldots &\ldots &\ldots &
}\]
\begin{center}
\refstepcounter{diagramm}
{\it Figure }\arabic{diagramm}: Tree labelled by the elements of $PSL_2(\ZZ)$\\[.5cm]
\label{trees}
\end{center}
\end{minipage}
}

We will first give the algorithm and then proof that the lists $\gen$ and $\rep$ that
are calculated
are what they should be.\\

\framebox{
\begin{minipage}{\linewidth}
{\bf Algorithm for Calculating $\boldmath{\bar{\Gamma}(O)}$:}\\[.2cm]
{\em Given:} Origami $O$.\\[.2cm]
Let $\rep$ and $\gen$ be empty lists.\\
Add $\bar{I}$ to $\rep$. $\bar{A} := \bar{I}$.\\
Loop:\\
$B := A\cdot T$,
$C := A\cdot S$\\
//Check whether $\bar{B}$ is already represented by $\rep$ and add, if there occur some, \\
//the found generators to $\gen$:\\
\hspace*{0.5cm}
\begin{minipage}{11cm}
For each $\bar{D}$ in $\rep$, check whether $B\cdot D^{-1}$ is in $\Gamma(O)$ or $-B\cdot D^{-1}$
is in $\Gamma(O)$. If so, add $\bar{B}\cdot\bar{D}^{-1}$ to $\gen$.\\
If none is found, add $\bar{B}$ to $\rep$.
\end{minipage}\\[.2cm]
Do the same for $C$ instead of $B$.\\
If there exists a successor of $\bar{A}$ in $\rep$, let $\bar{A}$ be now this successor and go to the beginning of the loop.\\
If not, finish the loop.\\[.2cm]
{\em Result:} $\gen$: list of generators of $\bar{\veech}(O)$, $\rep$: list of coset representatives
in $\pslzwei(\ZZ)$.
\end{minipage}
}
\newpage
\begin{rem}\hspace*{\fill}\\[-.5cm]
\begin{enumerate}
\item \label{differentcosets}
Any two elements of $\rep$ belong to different cosets.
\item \label{finite}
The algorithm stops after finitely many steps.
\item \label{rep}
In the end each coset is represented by a member of $\emrep$.
\item \label{gen}
In the end $\veechquer(O)$ is generated by the elements of $\emgen$.
\end{enumerate}
\end{rem}

\begin{proof}\hspace*{\fill}\\
\underline{\bf \ref{differentcosets}.:}
The statement follows by induction. It is true in the beginning, since $\rep$ contains only $\bar{I}$. After passing through the loop  it is still true, since $\bar{B}$ (respectively $\bar{C}$) is only added if $\bar{B}\cdot \bar{D}^{-1}$ (resp.\ $\bar{C}\cdot \bar{D}^{-1}$) is not in $\veechquer(O)$ for all $\bar{D}$ in $\rep$.\\[.2cm]
\underline{\bf \ref{finite}.:}
Follows from \ref{differentcosets}, since $\veechquer(O)$ has finite index in $\pslzwei(\ZZ)$ (\cite{gj}, Thm.\ 5.5).\\[.2cm]
\underline{\bf \ref{rep}.:} Let $\bar{A}$ be an arbitrary element
of $\pslzwei(\ZZ)$. There is at least one vertex in the tree that
is labelled by $\bar{A}$. Denote the vertices by $v_0$, $v_1$,
$v_2$, $\ldots$ as in Figure \ref{trees} and their
labels by $\bar{A}_0$, $\bar{A}_1$, $\bar{A}_2$, $\ldots$, respectively. \\
We do induction by the numeration $n$ of the vertices:\\
$\bar{A}_0 = \bar{I}$ is in $\rep$.
Suppose for a certain $n \in \NN$ all $\bar{A}_k$ with $k \leq n$ are represented by $\rep$.\\
If $A_{n+1}$ is not itself in $\rep$ then consider the path $\ww$ from $v_0$ to $v_{n+1}$
and let $v_j$ be the first vertex on $\ww$ that is not in $\rep$. Hence, its predecessor is
in $\rep$ and $\bar{A}_j$ was checked but not added. Thus, there is some
$\bar{A}_l$ ($l < j$) in $\rep$ such that $\bar{A}_j\cdot \bar{A}_l^{-1}$ is in
$\bar{\veech}(O)$, i.e. $\bar{A}_j \sim \bar{A}_l$.\\
Let $\wwdach$ be the path from $v_j$ to $v_{n+1}$ and $\bar{D}$ the product of the labels
of the edges on $\wwdach$. Then $\bar{A}_{n+1} = \bar{A}_{j}\cdot\bar{D}$.\\
Walking 'the same path' as $\wwdach$  starting at $v_l$ (i.e. a path described by the same
sequence of $\bar{S}$ and $\bar{T}$) leads to some vertex $v_m$ with $m < n+1$ and label
$\bar{A}_m = \bar{A}_l \cdot \bar{D}$. \\
We have $\bar{A}_{n+1} = \bar{A}_j\cdot\bar{D} \sim\bar{A}_l\cdot\bar{D} = \bar{A}_m$ and
by the assumption $\bar{A}_m$ is represented by $\rep$, hence also
$\bar{A}_{n+1}$ is.\\[.2cm]
\underline{\bf \ref{gen}.:}
Let $G$ be the group generated by the elements of $\gen$. We have by construction of the
list $\gen$ that $G \subseteq \veechquer(O)$.\\
We show again by induction that each label $\bar{A}_n$ in the tree that is in
$\veechquer(O)$ is also in $G$. This is true for $n = 0$. Suppose it is true for all
$k \leq n$ with a certain  $n\in\NN$.\\
If $\bar{A}_{n+1}$ is in $\veechquer(O)$, we proceed as in (\ref{rep})\ and find some $\bar{A}_j$, $\bar{A}_l$, $\bar{A}_m$ and $\bar{D}$ ($j,l,m < n+1$) such that $\bar{A}_j$ and $\bar{A}_l$ are in the same coset, $\bar{A}_j\cdot \bar{A}_l^{-1}$ is in the list $\gen$ (hence, $\bar{A}_j\cdot \bar{A}_l^{-1} \in G$), $\bar{A}_{n+1} = \bar{A}_j\cdot \bar{D}$ and $\bar{A}_m = \bar{A}_l\cdot \bar{D}$. $\bar{A}_m$ is in the same coset as $\bar{A}_{n+1}$, thus it is an element of $\veechquer(O)$. By the assumption $\bar{A}_m$ is then also in $G$. Hence, we have:
\[\bar{A}_{n+1} = \bar{A}_j\cdot \bar{A}_l^{-1}\cdot \bar{A}_l\cdot \bar{D} = (\bar{A}_j\cdot \bar{A}_l^{-1})\cdot \bar{A}_m\in G.\]
\end{proof}

Now - knowing $\bar{\veech}(O)$ -, it is easy to determine $\Gamma(O)$. We just have to distinguish the two cases,
whether $-I$ is in $\Gamma(O)$ or not.\\

\framebox{
\begin{minipage}{\linewidth}
{\bf Algorithm for Calculation of} $\boldmath{\veech(O)}$:\\[.2cm]
{\em Given:} Origami $O$.\\[.2cm]
Calculate $\gen$ and $\rep$.\\
Let $\gen'$ and $\rep'$ be empty lists.\\
Check, whether $-I \in \veech(O)$.\\
If yes:
\begin{minipage}[t]{9cm}
For each $\bar{A} \in \gen$ add $A$ to $\gen'$. Add $-I$ to $\gen'$.\\
For each $\bar{A} \in \rep$ add $A$ to $\rep'$.
\end{minipage} \\
If no:\phantom{s}
\begin{minipage}[t]{9cm}
For each $\bar{A} \in \gen$, check whether $A \in \veech(O)$.\\
If it is, add $A$ to $\gen'$; if it is not, add $-A$ to $\gen'$.\\
For each $\bar{A} \in \rep$ add $A$ and $-A$ to $\rep'$.
\end{minipage}\\[.2cm]
{\em Result}:
\begin{minipage}[t]{11cm}$
\gen'$: list of generators of $\veech(O)$,\\
$\rep'$: list of right coset representatives of $\veech(O)$ in $\slzwei(\ZZ)$.
\end{minipage}
\end{minipage}
}

\begin{ex} (for $O = L(2,3)$)\label{exgenerators}\\
1) Result of calculating $\bar{\veech}(O)$:\\[.2cm]
$\gen$:
\begin{eqnarray*}
\begin{pmatrix}1&3\\0&1\end{pmatrix} = \bar{T}^3,
\begin{pmatrix}-1&3\\-2&5 \end{pmatrix}= \bar{T}\bar{S}\bar{T}^2\bar{S}\bar{T}^{-1}\bar{T}^{-1},
\begin{pmatrix}1&0\\2&1 \end{pmatrix}= \bar{T}\bar{S}\bar{T}\bar{S}\bar{T}^{-1}\bar{S},
\begin{pmatrix}3&-5\\2&-3 \end{pmatrix}= \bar{T}^2\bar{S}\bar{T}\bar{S}\bar{T}^{-1}\bar{S}^{-1}\bar{T}^{-2}
\end{eqnarray*}
is a list of generators of $\bar{\veech}(O)$.\\
$\rep:$
\[
\bar{I},\, \bar{T},\, \bar{S},\, \bar{T}^2, \,\bar{T}\bar{S},\,
\bar{S}\bar{T},\, \bar{T}^2\bar{S},\,
\bar{T}\bar{S}\bar{T},\, \bar{T}^2\bar{S}\bar{T}
\]
is a system of coset representatives of $\bar{\veech}(O)$ in $\slzwei(\ZZ)$.\\
(The algorithm produces more generators (compare example \ref{exfundamental}). We eliminated
redundant ones.)\\

2) Result of calculating $\veech(O)$: ($-I \in \Gamma(O)$)
\begin{eqnarray*}
\gen' &=& \gen \cup \{-I\}. \\
\rep':&=&
I,\, T,\, S,\, T^2, \,TS,\,
 ST,\, T^2S,\,
 TST,\, T^2ST
\end{eqnarray*}

Hence, $\Gamma(O)$ is a subgroup of index $9$ in $\slzwei(\ZZ)$.
\end{ex}

\subsection{Geometrical type of $\boldmath{\HH/\bar{\veech}(O)}$}
\label{geometry}\hspace*{\fill}\\

The group $\bar{\veech}(O)$ is a subgroup of $\pslzwei(\ZZ)$ and
of finite index (\cite[Thm. 5.5]{gj}), thus it operates as
Fuchsian group (via M\"obius transformations) on $\HH$ and  $V :=
\HH/\bar{\Gamma}(O)$ is an affine algebraic curve. It is
defined over $\bar{\QQ}$ by the Theorem of Belyi:
We have a covering from $\HH/\bar{\veech}(O)$ to
$\HH/\pslzwei(\ZZ)\cong \AAA^1(\CC) = \PP^1(\CC)-\{\infty\}$
ramified at most over the images of $i$ and
$\rho = \frac12+(\frac12\sqrt{3})i$. Thus, by Belyi's theorem
the projective curve $\overline{\HH/\bar{\veech}(O)}$ and hence also $C$ is defined over $\Qquer$.\\
We want to determine the genus and the number of points
at infinity of the curve $\HH/\bar{\Gamma}(O)$.\\[.2cm]
Let $\Delta := \Delta(P_0,P_1,P_{\infty})$ be the
standard fundamental domain of $\slzwei(\ZZ)$, i.e. the hyperbolic
pseudo-triangle with vertices $P_0 := -\frac12 + \frac{\sqrt{3}}2i$, $P_1 := \frac12 +
\frac{\sqrt{3}}{2}i$ and
$P_{\infty}:= i\infty$.\\
We denote by $\bar{A}$ also the M\"obius transformation
defined by the matrix $A$. Then $\bar{T}$ and $\bar{S}$ (as M\"obius transformations)
send $P_0P_{\infty}$ to $P_1P_{\infty}$, respectively
$P_0P_1$  to itself (fixing $i$).\\
Let $\rep = \{\bar{A}_1,\ldots,\bar{A}_k\}$ be the system of right coset representatives we calculated in section
\ref{final}. Then
\[F := \bigcup_{i=1}^{k}\bar{A}_i(\Delta)\]
is a simply connected fundamental domain of $\bar{\veech}(O)$.
The list $\gen$ of generators defines how to
glue the edges of $F$ to obtain $\HH/\bar{\veech}(O)$.
This way, we get a triangulation of $\HH/\bar{\veech}(O)$
(compare Figure \ref{fundamentaldomain}).
We calculate the numbers $t$, $e$, $v$ of the triangles,
the edges and the vertices of this triangulation
as described in the following algorithm. Furthermore,
the vertices defined by translates of $P_{\infty}$
are exactly the cusps of $\HH/\bar{\veech}(O)$. We
denote their number by $\hat{v}$. Thus (using the formula of Euler for calculating the genus)
we get the following result.
\begin{rem}
Let $t$, $e$, $v$ and $\hat{v}$ be the numbers of
triangles, edges, vertices and marked vertices as calculated in the
following algorithm. Then $\HH/\bar{\veech}(O)$ is an affine curve of
genus $g = \frac{2-(v-e+t)}{2}$ with $\hat{v}$ cusps.
\end{rem}

\framebox{
\begin{minipage}{\linewidth}
{\bf Algorithm determining the geometrical type of $\boldmath{\HH/\bar{\Gamma}(O)}$:}\\[.2cm]
Generate a list of triangles  $L := \{\bar{A}_1(\Delta), \ldots, \bar{A}_k(\Delta)\}$.\\
In the triangle $\bar{A}_i(\Delta)$ we call $\bar{A}_i(P_0)\bar{A}_i(P_1)$ (the
image of the edge $P_0P_1$) 'the $S$-edge'.
Similarly, we call $\bar{A}_i(P_1)\bar{A}_i(P_{\infty})$ 'the $T$-edge'
and $\bar{A}_i(P_0)\bar{A}_i(P_{\infty})$ 'the $T^{-1}$-edge'.\\
For each $i,j \in \{1,\ldots,k\}$ identify
\begin{itemize}
\item
the $T$-edge of $\bar{A}_j(\Delta)$ with the $T^{-1}$-edge of $\bar{A}_i(\Delta)$, if $\bar{A}_i \sim \bar{A}_j
\cdot \bar{T}$,
i.e. if $(\bar{A}_j\bar{T})\bar{A}_i^{-1} \in \bar{\veech}(O)$,
\item
the $T^{-1}$-edge of $\bar{A}_j(\Delta)$ with the $T$-edge of $\bar{A}_i(\Delta)$, if $\bar{A}_i \sim
\bar{A}_j\cdot \bar{T}^{-1}$
and
\item
the $S$-edge of $\bar{A}_j(\Delta)$ with the $S$-edge of $\bar{A}_i(\Delta)$, if $\bar{A}_i \sim
\bar{A}_j\cdot \bar{S}$.\\
If an $S$-edge of some triangle $\bar{A}_j(\Delta)$ is identified with
itself (i.e. $i = j$) create an additional triangle: Add a vertex in the middle of this $S$-edge and
add an edge from this new vertex to the opposite vertex in the
triangle $\bar{A}_j(\Delta)$.  (Compare triangle $\bar{T^2}\bar{S}\bar{T}$ in
Figure \ref{fundamentaldomain}). This is done
to get in the end a triangulation of the surface.
\end{itemize}
$t :=$ number of triangles. \quad
$e :=$ number of edges.\\
$v :=$ number of vertices, $\hat{v} := $ number of vertices that are endpoints of $T$-edges.\\
$g := \frac{2-(v-e+t)}{2}$.\\[.2cm]
{\em Result:}
\begin{minipage}[t]{11cm}
$g$ : genus of $\HH/\bar{\veech}(O)$ \quad
$\hat{v}$: number of vertices at infinity of $\HH/\bar{\veech}(O)$.
\end{minipage}
\end{minipage}
}

\begin{ex} \label{exfundamental} (for $O = L(2,3)$)\\ \nopagebreak

$\rep$: $\bar{I}, \bar{T}, \bar{T}^2, \bar{T}^2\bar{S}, \bar{T}^2\bar{S}\bar{T}, \bar{T}\bar{S}, \bar{T}\bar{S}\bar{T}, 
\bar{S}, \bar{S}\bar{T}$.\\
$\gen$:
\begin{minipage}[t]{11cm}
$a := \bar{T}^3, b := \bar{S}\bar{T}\bar{S}\bar{T}^{-1}\bar{S}\bar{T}^{-1}, c := \bar{S}\bar{T}^2\bar{S}, 
d := \bar{T}\bar{S}\bar{T}^2\bar{S}\bar{T}^{-2},$\\ $e := \bar{T}\bar{S}\bar{T}^{-2}\bar{S}\bar{T}^{-2},$
$f := \bar{T}^2\bar{S}\bar{T}\bar{S}\bar{T}^{-1}\bar{S}\bar{T}^{-2}$
\end{minipage}\\
\includegraphics[scale = 0.4]{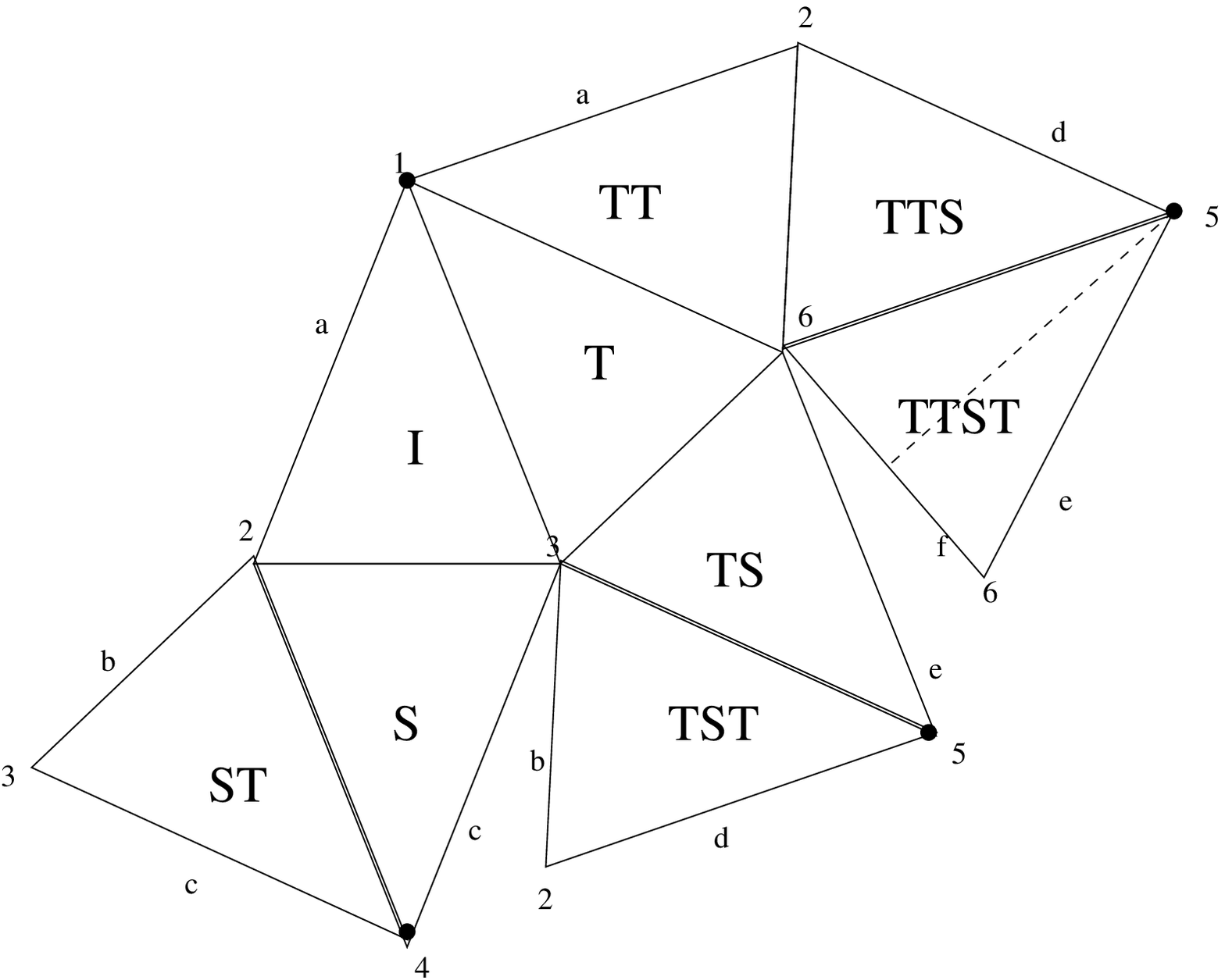}
\begin{center}
\refstepcounter{diagramm}
{\it Figure }\arabic{diagramm}: Fundamental domain of $\bar{\Gamma}(L(2,3))$.\\[.5cm]
\label{fundamentaldomain}
\end{center}
Edges with the same letters are glued. In triangle $\bar{T}^2\bar{S}\bar{T}(\Delta)$ an edge and a vertex were added, since the
'S-edge' is glued to itself. Vertices with same numbers are identified. Vertices at infinity are marked by
a filled circle.\\
Thus, $t = 9+1, e = 14 + 1, v = 6 + 1, \hat{v} = 3$.\\

\framebox{
\begin{minipage}{11cm}
{\bf Result:}  $g = 0, \hat{v} = 3$. Hence,
\[\HH/\bar{\Gamma}(L(2,3)) \cong \PP^1 - \{0,1,\infty\}.\]
\end{minipage}
}
\end{ex}

\begin{prop}\label{cong}
$\veech(L(2,3))$ is not a congruence subgroup of $\emslzwei(\ZZ)$ .
\end{prop}

\begin{proof}
Recall from Example \ref{exgenerators} that
\[\veech(L(2,3)) =
   <\begin{pmatrix}1&3\\0&1\end{pmatrix}, 
    \begin{pmatrix}1&0\\2&1\end{pmatrix}, \begin{pmatrix}-1&3\\-2&5\end{pmatrix}, 
    \begin{pmatrix} 3&-5\\2&-3\end{pmatrix},
    \begin{pmatrix} -1&0\\0&-1\end{pmatrix}>.\]
$\HH/\bar{\veech}(L(2,3))$ has three cusps represented in Figure \ref{fundamentaldomain} by the vertices
1, 4 and 5. $T^{3}$, $ST^2S^{-1}$ and $TST^4S^{-1}T^{-1}$ are parabolic elements that correspond to them
respectively and the amplitudes are 3, 2 and 4. Hence, the
level $m$ of  $\veech(L(2,3))$ is lcm$(3,2,4) = 12$ (using notations of \cite{wo}).\\
Suppose that $\veech(L(2,3))$ is a congruence subgroup. By Theorem 2 in \cite{wo} we would have:
\begin{equation}
\Gamma(12) \subseteq \veech(L(2,3)). \label{congruence}
\end{equation}
Let $p : \pslzwei(\ZZ) \rightarrow \pslzwei(\ZZ/3\ZZ)$ be the natural projection.
Then we have
\[p(\bar{\veech}(L(2,3))) = <\begin{pmatrix}\bar{1}&\bar{0}\\ \bar{2}&\bar{1}\end{pmatrix},
                       \begin{pmatrix}\bar{0}&\bar{1}\\\bar{2}&\bar{0}\end{pmatrix}>
                    = \pslzwei(\ZZ/3\ZZ).\]
Hence Diagram \ref{noncongruencel23} is commutative with
$N := \bar{\veech}(L(2,3)) \cap \bar{\Gamma}(3)$.\\
\[ \xymatrix{
   1 \ar[r]& \bar{\Gamma}(3) \ar[r]& \pslzwei(\ZZ)  \ar[r] & \pslzwei(\ZZ/3\ZZ)  \ar[r] & 1\\
   1 \ar[r]& N \ar[r] \ar@{^{(}->}[u] &\bar{\veech}(L(2,3))  \ar[r] \ar@{^{(}->}[u] &
     \pslzwei(\ZZ/3\ZZ)  \ar[r] \ar@{=}[u]& 1
   }
\]
\begin{center}\refstepcounter{diagramm}{\it Diagram }\arabic{diagramm}\label{noncongruencel23}
\end{center}
Since the index $[\pslzwei(\ZZ):\bar{\veech}(L(2,3))]$ of $\bar{\veech}(L(2,3))$
in $\pslzwei(\ZZ)$ is $9$ it follows from Diagram \ref{noncongruencel23} that
$[\bar{\Gamma}(3): N] = 9$.\\
By (\ref{congruence}) we have: $\bar{\Gamma}(12) \subseteq N \subseteq \bar{\Gamma}(3)$. But
$[\bar{\Gamma}(3):\bar{\Gamma}(12)] =  2^4\cdot3$
(using \cite{sh}, (1.6.2)). Thus $[\bar{\Gamma}(3): N] = 9$ would have to be a factor of $2^4\cdot3$. 
Contradiction!

\end{proof}

\section{Some examples}

\subsection{"Trivial Origamis":}\label{trivial}\hspace*{\fill}\\

\setlength{\unitlength}{0.5cm}
\begin{minipage}{4.5cm}
\hspace*{1.5cm}
\begin{picture}(4.3,3.3)
\put(-2.5,1.2){$O =$}
\put(0,0){\framebox(1,1){}}
\put(1,0){\framebox(1,1){}}
\put(2,0){\framebox(1,1){}}
\put(3,0){\framebox(1,1){}}
\put(0,1){\framebox(1,1){}}
\put(1,1){\framebox(1,1){}}
\put(2,1){\framebox(1,1){}}
\put(3,1){\framebox(1,1){}}
\put(0,2){\framebox(1,1){}}
\put(1,2){\framebox(1,1){}}
\put(2,2){\framebox(1,1){}}
\put(3,2){\framebox(1,1){}}
\put(0.2,-0.4){$a_1$}
\put(1.6,-0.4){$\ldots$}
\put(3.2,-0.4){$a_n$}
\put(0.2,3.2){$a_1$}
\put(1.6,3.2){$\ldots$}
\put(3.2,3.2){$a_n$}
\put(-0.9,0.2){$b_1$}
\put(-0.7,1.2){$\vdots$}
\put(-0.9,2.2){$b_m$}
\put(4.2,0.2){$b_1$}
\put(4.4,1.2){$\vdots$}
\put(4.2,2.2){$b_m$}
\end{picture}
\end{minipage}
\phantom{t}
\begin{minipage}{4cm}
\begin{eqnarray*}\Gamma(O) = \{\begin{pmatrix}a&b\\c&d\end{pmatrix}
\in \slzwei(\ZZ)|b\equiv 0 \mbox{ mod }n', c\equiv 0 \mbox{ mod }m'\}\\
\mbox{ where } t := \mbox{ gcd}(m,n), n' := n/t, m' := m/t
\end{eqnarray*}\\
\end{minipage}

\subsection{"$L$-Sequence":}\label{lserie}\hspace*{\fill}\\

\hspace*{2.5cm}
\begin{minipage}{3cm}
\setlength{\unitlength}{0.5cm}
\begin{picture}(5,3.3)
\put(-4,1.2){$L(n,m) =$}
\put(0,0){\framebox(1,1){}}
\put(1,0){\framebox(1,1){}}
\put(2,0){\framebox(1,1){}}
\put(3,0){\framebox(1,1){}}
\put(0,1){\framebox(1,1){}}
\put(0,2){\framebox(1,1){}}
\put(0.2,-0.4){$a_1$}
\put(1.1,-0.4){$a_2$}
\put(1.5,-0.4){$\ldots$}
\put(3.2,-0.4){$a_n$}
\put(0.2,3.2){$a_1$}
\put(1.5,1.2){$\ldots$}
\put(3.2,1.2){$a_n$}
\put(-0.9,0.2){$b_1$}
\put(-0.7,1.5){$\vdots$}
\put(-0.9,2.4){$b_m$}
\put(4.2,0.2){$b_1$}
\put(1.4,1.5){$\vdots$}
\put(1.2,2.4){$b_m$}
\end{picture}
\end{minipage}
\phantom{t}
\begin{minipage}{4cm}
\begin{tabular}{l|l|l|l}
Origami&Index&Genus&$\sharp$ Cusps\\
\hline
$L(2,2)$&3&0&2\\
$L(2,3)$&9&0&3\\
$L(2,4)$&18&0&5\\
$L(2,5)$&36&0&8\\
$L(2,6)$&54&0&10\\
$L(2,7)$&108&1&17\\
$L(3,3)$&9&0&3\\
$L(4,4)$&54&0&10\\
\end{tabular}
\end{minipage}\\

\subsection{"Cross - Sequence":}\label{kreuzserie}\hspace*{\fill}\\
\hspace*{2cm}
\begin{minipage}{4.5cm}
\begin{center}
\setlength{\unitlength}{.6cm}
\begin{minipage}{5cm}
\begin{picture}(4,2.3)
\put(-2.5,0.3){$O_{2k} =$}
\put(0,0){\framebox(1,1){1}}
\put(1,0){\framebox(1,1){}}
\put(2,0){\framebox(1,1){}}
\put(3,0){\framebox(1,1){}}
\put(4,0){\framebox(1,1){}}
\put(5,0){\framebox(1,1){2k}}
\put(0.2,-0.4){$a_1$}
\put(1.1,-0.4){$a_2$}
\put(2.8,-0.4){$\ldots$}
\put(3.8,-0.4){$a_{2k-1}$}
\put(5.4,-0.4){$a_{2k}$}

\put(0.2,1.2){$a_2$}
\put(1.1,1.2){$a_1$}
\put(2.8,1.2){$\ldots$}
\put(4,1.2){$a_{2k}$}
\put(5.1,1.2){$a_{2k-1}$}

\put(-0.6,0.5){$a_0$}
\put(6.2,0.5){$a_0$}

\end{picture}
\end{minipage}\\[1.5cm]
\end{center}
\end{minipage}
\begin{minipage}{4cm}
\begin{center}
\begin{tabular}{l|l|l|l}
Origami&Index&Genus&$\sharp$ Cusps\\
\hline
$O_2$&3&0&2\\
$O_4$&6&0&3\\
$O_6$&12&0&4\\
$O_8$&24&0&6\\
$O_{10}$&36&0&8\\
$O_{12}$&48&0&10\\
$O_{14}$&72&1&12\\
$O_{16}$&96&2&14\\
\end{tabular}\\[1cm]
\end{center}
\end{minipage}\\

\subsection{Remarks:}\hspace*{\fill}\\

As in Example \ref{l23} edges labelled with same letters are
glued. The tables in \ref{lserie} and \ref{kreuzserie} itemize for
an origami $O$ respectively the index of the projective Veech
group $\bar{\Gamma}(O)$
in $\pslzwei(\ZZ)$ and the genus and number of cusps of $\HH/\bar{\veech}(O)$.\\
For the example in \ref{trivial}, $\veech(O)$ can be determined using Proposition
\ref{mainprop}.\\
The sequence in \ref{lserie} was introduced to me by Pierre 
Lochak. The Veech group
e.g.\ of $L(2,2)$ is given also  in \cite{mm}. This sequence
is also studied in detail in \cite{li} and e.g.\ estimates for the growth 
of the genus
and the number of cusps are obtained.
The Veech groups in this sequence are in general not congruence subgroups of
$\slzwei(\ZZ)$ (see Proposition \ref{cong}).\\
On the contrary one can show - again using Proposition \ref{mainprop} - that the Veech groups
$\veech(O_{2k})$ in \ref{kreuzserie} are congruence
subgroups for all $k \in \NN$. Furthermore the genus of the curve
$\HH/\veech(O_{2k})$ is not bounded.\footnote{Details
    will be published elsewhere}\\

Only a few general statements about Veech groups of origamis are known yet.
There seems to be no obvious relation between the index $d$ of the origami
$O = (p: X \rightarrow \torus)$ and the index of its Veech group.
In particular, it follows from Proposition \ref{mainprop} that
each characteristic subgroup of $F_2$ defines an origami
with Veech group $\slzwei(\ZZ)$. (The smallest, nontrivial example
(calculated by Frank Herrlich) is defined by a covering $p: X
\rightarrow  \torus$ of degree $108$.)
Hence, there is a cofinal system of origamis having the full
group $\slzwei(\ZZ)$ as Veech group.

\end{document}